\nonstopmode
\input amstex
\input amsppt.sty   
\hsize 12.3 cm 
\vsize 47pc
\def\nmb#1#2{#2}         
\def\cit#1#2{\ifx#1!\cite{#2}\else#2\fi} 

\define\De{\Delta}

\define\row#1#2#3{#1_{#2},\ldots,#1_{#3}}

\def\today{\ifcase\month\or
 January\or February\or March\or April\or May\or June\or
 July\or August\or September\or October\or November\or December\fi
 \space\number\day, \number\year}
\message{")"}

\topmatter
\title Choosing roots of polynomials smoothly, II
\endtitle
\author 
Andreas Kriegl, Mark Losik, and Peter W\. Michor  \endauthor
\leftheadtext{A\. Kriegl, M\. Losik, P\.W\. Michor}
\address
A\. Kriegl: Institut f\"ur Mathematik, Universit\"at Wien,
Strudlhofgasse 4, A-1090 Wien, Austria
\endaddress
\email Andreas.Kriegl\@univie.ac.at \endemail
\address M\. Losik: Saratov State University, ul\. Astrakhanskaya, 83, 
410026 
Saratov, Russia \endaddress
\email losik\@info.sgu.ru\endemail
\address
P\. W\. Michor: Institut f\"ur Mathematik, Universit\"at Wien,
Strudlhofgasse 4, A-1090 Wien, Austria; {\it and:}
Erwin Schr\"odinger Institut f\"ur Mathematische Physik,
Boltzmanngasse 9, A-1090 Wien, Austria
\endaddress
\email Peter.Michor\@esi.ac.at \endemail

\thanks{M.L. and P.W.M. were supported
     by `Fonds zur F\"orderung der
     wissenschaftlichen For\-schung,
     Projekt P~14195~MAT'.}
\endthanks
\abstract 
        We show that the roots of any smooth curve of polynomials
        with real roots only 
        can be parametrized twice differentiable (but not better).
\endabstract
\endtopmatter
\document

In \cit!{1} we claimed that there exists a smooth curve of polynomials
of degree 3 for which no $C^1$-parametrization of the roots exists.
Unfortunately there was an error in the calculation of $b_3$ and we have
been informed by Jacques Chaumat and Anne-Marie Chollet in June 2001 
about that and the related papers \cit!{2}, \cit!{5}.
 
We are now going to repair this mistake and improve at the same time 
the results of \cit!{2}.
The smoothness assumptions in the following theorem
are certainly not the best possible but in fact we are mainly 
interested in the case of smooth coefficients.

The conclusion of the theorem is the best possible, since even for
the characteristic polynomial of a smooth curve of symmetric matrices
there needn't be a differentiable parametrization of the roots 
with locally H\"olderian derivative as the
first example in \cit!{3} shows.

Let $P$ be a curve defined on some subset $T\subseteq \Bbb R$ of monic 
polynomials $P(t)$ of degree $n\geq 1$ with real roots only.
A parametrization of some class of the roots of $P$ is a curve 
$x:T\to\Bbb R^n$ of that class such that for each $t\in T$ the values 
$x_1(t),\dots,x_n(t)$ are the roots of $P(t)$ with correct multiplicity.

\proclaim{Theorem}
Consider a continuous curve of polynomials
$$ 
P(t)(x)=x^n-a_1(t)x^{n-1}+\dots+(-1)^na_n(t), \quad t\in \Bbb R,
$$ 
with all roots real.
Then there is a continuous parametrization 
$x=(x_1,\dots,x_n):\Bbb R\to\Bbb R^n$ of the roots of $P$.
Moreover:
\roster
\item"(\nmb:{1})" \cit!{2},~Theorem~1 and Theorem~2.
       If all coefficients $a_i$ are of class 
       $C^{n}$ then the parametrization $x:\Bbb R\to\Bbb R^n$ 
       may be chosen differentiable with locally bounded derivative. 
\item"(\nmb:{2})" If all $a_i$ are of class $C^{2n}$ then any 
       differentiable parametrization $x:\Bbb R\to\Bbb R^n$ is 
       actually $C^1$. 
\item"(\nmb:{3})" If all $a_i$ are of class $C^{3n}$ then the 
       parametrization $x:\Bbb R\to\Bbb R^n$ may be chosen twice 
       differentiable. 
\endroster
\endproclaim

\demo{Proof} The parameterization by order
$x_1(t)\le\dots\le x_n(t)$ is continuous, see e.g\. 
\cit!{1},~4.1.
We prove \therosteritem{\nmb|{2}} and 
\therosteritem{\nmb|{3}}, and
we use the proof of theorem 4.3 in \cit!{1}.
First we replace $x$ by $x+\frac1na_1(t)$, and consequently assume 
without loss that $a_1=0$.

As noted in the proof of 4.3 in \cit!{1}
the multiplicity lemma \cit!{1},~\nmb!{3.7} remains true in 
the $C^m$-case for $m\ge n$ in the following sense, with the same proof:
\newline
{\sl If $a_1=0$ then the following two conditions are equivalent
\roster
\item $a_k(t)=t^k a_{k,k}(t)$ for a $C^{m-k}$-function $a_{k,k}$, 
       for all $2\le k\le n$.
\item $a_2(t)=t^2 a_{2,2}(t)$ for a $C^{m-2}$-function $a_{2,2}$.
\endroster}

\medskip
\noindent {\it Proof of \therosteritem{\nmb|{2}}.}
Let all $a_i$ be $C^{2n}$. 

Then we choose a fixed $t$, say $t=0$. 

If $a_2(0)=0$ then it vanishes of second order at 0: if it vanishes only 
of first order then $\tilde\De_2(P(t))=-2na_2(t)$ (see \cit!{1},~3.1)
would change sign at $t=0$, contrary to the assumption that all roots 
of $P(t)$ are real, by \cit!{1},~\nmb!{3.2}.
Thus $a_2(t)=t^2a_{2,2}(t)$, so by the variant of the multiplicity 
lemma described above we have $a_k(t)=t^k a_{k,k}(t)$ for 
$C^n$-functions $a_{k,k}$, for $2\le k\le n$. 
We consider the following $C^n$-curve of polynomials 
$$ 
P^1(t)(z)=z^n+a_{2,2}(t)z^{n-2}-a_{3,3}(t)z^{n-3}\dots 
+(-1)^na_{n,n}(t).
$$
Then $P(t)(tz)=t^n\,P^1(t)(z)$ and hence $z\mapsto t\,z=x$ gives for 
$t\ne 0$ a bijective correspondance between the roots $z$ of $P^1(t)$ 
and the roots $x$ of $P(t)$ with correct multiplicities.
Moreover parametrizations $z$ which are continuous at $t=0$
correspond to parametrizations $x$ which are differentiable at $t=0$.
By \therosteritem{\nmb|{1}} we may choose the parametrization 
$z=(z_1,\dots,z_n)$ differentiable with locally bounded derivative.
Then the corresponding parametrization $t\mapsto x(t):=t\,z(t)$ is 
differentiable with derivative $x'(t)=t\,z'(t)+z(t)$ which is 
continuous at $t=0$ with $x'(0)=z(0)$.

If $a_2(0)\ne 0$ then we use the splitting lemma \cit!{1},~\nmb!{3.4} 
for the $C^{2n}$-case: 
We may factor $P(t)=P_1(t)\dots P_k(t)$ for $t$ in a neighborhood 
of $0$ and some $k>1$ where the $P_i$ have again 
$C^{2n}$-coefficients and where each $P_i(0)$ has all roots equal to, 
say, $c_i$, and where the $c_i$ are distinct.
By the argument above applied to each $P_i$ separately, there is a 
differentiable parametrization $x=(x_1,\dots,x_n)$ of roots whose  
derivative $x'$ is continuous at $t=0$.
Moreover, if $P_i(0)(x_j(0))=0$ then $x_j'(0)$ is a root of the 
polynomial $P_i^1(0)$ which depends only on $P_i$. 
We shall use this for arbitrary $t$ below.

{\it Claim.
Any differentiable parametrization $y=(y_1,\dots,y_n)$ of the roots of $P$
has $y'$ continuous at $t=0$:}
Let $i\in\{1,\dots,n\}$.
For $t_m\to 0$ 
there are $k_m\in\{1,\dots,n\}$ such that $y_i(t_m)=x_{k_m}(t_m)$.
Choose a subsequence of the $t_m$ again denoted $t_m$ such that 
$y_i(t_m)=x_k(t_m)$ for some fixed $k$ and all $m$. 
By the argument above then we also have $y_i'(t_m)=x_{j_m}'(t_m)$ for 
some $j_m$ with $x_{j_m}(t_m)=x_k(t_m)=y_i(t_m)$. 
Passing again to a subsequence we find a fixed $j$ 
such that $y_i(t_m)=x_j(t_m)$ and $y_i'(t_m)=x_j'(t_m)$.
Then 
$$\gather
y_i(0) = \lim_m y_i(t_m) = \lim_m x_j(t_m) = x_j(0)\\
y_i'(0) = \lim_m \frac{y_i(t_m)-y_i(0)}{t_m} = 
\lim_m\frac{x_j(t_m)-x_j(0)}{t_m} = x_j'(0)  
\endgather$$
and so $y_i'(t_m)=x_j'(t_m)\to x_j'(0)=y_i'(0)$. 

Thus any differentiable parametrization of the roots of $P$ 
(which exists by \therosteritem{\nmb|{1}}) is indeed $C^1$, and 
\therosteritem{\nmb|{2}} is proved.

\medskip

\noindent {\it Proof of \therosteritem{\nmb|{3}}.}
Let all $a_i$ be $C^{3n}$. 
Remember that $a_1=0$.

\therosteritem{a}
Choose a fixed $t$, say $t=0$. 
If $a_2(0)=0$ then we consider again the polynomials $P^1(t)$, which 
now form a $C^{2n}$-curve. 
By \therosteritem{\nmb|{2}} its roots can be parametrized by a 
$C^1$-curve $t\mapsto z(t)=(z_1(t),\dots,z_n(t))$. 
The $x(t)=t\,z(t)$ are then again the roots of $P(t)$, now with 
continuous derivative $x'(t)=t\,z'(t)+z(t)$ which is differentiable  
at $t=0$ with $x''(0)= 2\,z'(0)$.

We show by induction on $n$ that for fixed open intervalls 
$I\subseteq\Bbb R$ there exists a twice differentiable 
parametrization $y$ of the roots of $P$ on $I$. 

Let $t_0\in I$ be such that $a_2(t_0)\ne 0$.
By the splitting lemma \cit!{1},~\nmb!{3.4} for the $C^{3n}$-case 
we may factor $P(t)=P_1(t)\dots P_k(t)$ for some $k>1$ and all 
$t$ in a neighborhood $I_1\subseteq I$ of $t_0$ where
the $P_i(t)$ have again $C^{3n}$-coefficients and where each 
$P_i(t_0)$ has all roots equal to, say, $c_i$, and where the $c_i$ 
are distinct.
By induction there is on $I_1$ a twice differentiable 
parametriziation of the roots of each $P_i$. 
Note that for $n=1$ the root equals the (single) coefficient.

Let now $a_2(t)\ne 0$ for all $t\in I$.
We consider twice differentiable parametrizations of the roots 
defined on open subintervalls $I_1\subseteq I$.
Obviously we may apply Zorn's lemma to obtain a twice differentiable 
parametrization on some maximal open subintervall $I_1$.
Suppose for contradiction that $I\supsetneq I_1$ and let the, say 
right, endpoint $t_0$ of $I_1$ belong to $I$.
Then there is a twice differentiable parametrization $y$ on $I_1$ and 
since $a_2(t_0)\ne 0$ a twice differentiable parametrization $x$ in a 
neighborhood of $t_0$. 
Let $t_m\nearrow t_0$.
For every $m$ there exists a permutation $\pi$ of $\{1,\dots,n\}$ 
such that $y_{\pi(i)}(t_m)=x_i(t_m)$ for all $i$. 
By passing to a subsequence, again denoted $t_m$, we may assume
that the permutation does not depend on $m$.
By passing again to a subsequence we may also assume
that $y_{\pi(i)}'(t_m)=x_i'(t_m)$ and then again for a subsequence 
that $y_{\pi(i)}''(t_m)=x_i''(t_m)$ for all $i$ and all $m$.
So we may paste $(y_{\pi(i)}(t))_i$ for $t<t_0$
with $x(t)$ for $t\geq t_0$ to obtain a twice differentiable 
parametrization on an intervall larger than $I_1$, a contradiction. 

Now we consider the closed set 
$E=\{t\in I\colon a_2(t)=0\}=\{t\in I\colon x_1(t)=\dots=x_n(t)\}$.  
Then $I\setminus E$ is open, thus a 
disjoint union of open intervals on which we have a twice differentiable
parametrization $x$ of the roots by the previous paragraph.

Consider next the set $E'$ of all accumulation points of $E$.
Then $I\setminus E'=(I\setminus E)\cup (E\setminus E')$ is again open 
and thus a disjoint union of open intervals, and for each point 
$t_0\in E\setminus E'$, i.e\. isolated point of
$E$, we have a twice differentiable local parametrization of roots 
$y_i(t)$ for $t\ne t_0$ (left and right of $t_0$), and we have a 
local $C^1$ parametrization $x_k(t)$ for $t$ near $t_0$ which is 
twice differentiable at $t_0$, by argument \therosteritem{a}. 
Clearly $y_i(t)\to x_1(t_0)=\dots=x_n(t_0)$ for $t\to t_0$. 

For $t_m\searrow t_0$, by passing to a subsequence, we may assume 
that $y_i'(t_m)=x_{\pi(i)}'(t_m)\to x_{\pi(i)}'(t_0)$.
Thus $y_i'(t)$ has at most $x_1'(t_0),\dots x_n'(t_0)$ as cluster 
points for $t\searrow t_0$.
Since $y_i'$ satisfies the intermediate value theorem, $y_i'(t)$ 
converges for $t\searrow t_0$, with limit $x_{\pi(i)}'(t_0)$, since 
it does so along a sequence $t_m$ as above. 
By renumbering the $y_i$ to the right of $t_0$ we may assume that 
$i=\pi(i)$. 
Similarly for the left side of $t_0$. 
Then $y_i'(t)\to x_i'(t_0)$ for $t\to t_0$, so $y_i$ is $C^1$ near 
$t_0$ and still twice differentiable off $t_0$.

In order to get twice differentiability at $t_0$ also, we consider 
again the situation at the beginning of the last paragraph.
Then we have
$$
\frac{y_i'(t_m) - y_i'(t_0)}{t_m-t_0} = \frac{x_{\pi(i)}'(t_m) - 
x_{\pi(i)}'(t_0)}{t_m-t_0} \to x_{\pi(i)}''(t_0) 
$$
so that 
$(y_i'(t)-y_i'(t_0))/(t-t_0)$ has at most 
$\{x_j''(t_0): x_j'(t_0)=y_i'(t_0)\}$ as cluster points for 
$t\searrow t_0$.
Since it satisfies the intermediate value theorem it 
converges for $t\searrow t_0$, with limit $x_{\pi(i)}''(t_0)$, since 
it does so along a sequence $t_m$ as just used.
Similarly for the left handed second derivative.
Thus we may renumber those $y_i$ for which the $y_i'(t_0)$ agree, to 
the right of $t_0$ in such a way that the (one sided) second 
derivatives agree. 
Then the (twice) renumbered $y_i$ are twice differentiable also at $t_0$.   

Thus we have a twice differentiable parametrization of roots on the open 
set $I\setminus E'$. 

Now let $t_0\in E'$, i.e\. an accumulation point of $E$. 
Let $F$ the set of all $t\in I$ where $x_1(t)=\dots=x_n(t)$ and 
$x_1'(t)=\dots=x_n'(t)$.  
Then $t_0\in F$ since each $x_i'(t_0)$ may be computed using only 
points in $E$.
Let $F'$ be the set of all accumulation points of $F$. 
Thus $E'\subseteq F= (F\setminus F')\cup F'\subseteq E$.

Let first $t_0\in F\setminus F'$, i.e\. an isolated point in $F$. 
Then again we have a local twice differentiable parametrization 
$t\mapsto y(t)$ of the roots for $t\ne t_0$ (left and right of 
$t_0$), since near $t_0$ there are only points in $I\setminus E'$.  
We still have a local $C^1$ parametrization $x$ near $t_0$ which is 
twice differentiable at $t_0$, by the argument above. 
As above we can find a twice differentiable parametrization $y$ of 
the roots on the open set $(I\setminus E')\cup(F\setminus F')$.

Finally, let $t_0\in F'$, i.e\. an accumulation point in $F$. 
We use again parameterizations $x$ near $t_0$, and $y$ as above.
Then all $x_i(t_0)$ agree, all $x_i'(t_0)$ agree, and even all 
$x_i''(t_0)$ agree.
We extend each $y_i$ from $(I\setminus E')\cup(F\setminus F')$ 
by these single function on $F'$ to the whole of 
$(I\setminus E')\cup (F\setminus F')\cup F'=(I\setminus E')\cup F=I$. 
We have to check that then each $y_i$ is 
twice differentiable at $t_0$. 
For $t_m\to t_0$ we have, by passing to a subsequences,
$$\align
y_i(t_m) = x_j(t_m) &\to x_j(t_0) = x_i(t_0) = y_i(t_0) 
\\ 
\frac{y_i(t_m) - y_i(t_0)}{t_m-t_0} = \frac{x_j(t_m) - 
x_j(t_0)}{t_m-t_0} &\to x_j'(t_0) =x_i'(t_0)
\\ 
\frac{y_i'(t_m) - y_i'(t_0)}{t_m-t_0} = \frac{x_j'(t_m) - 
x_j'(t_0)}{t_m-t_0} &\to x_j''(t_0) = x_i''(t_0) \qed
\endalign$$
\enddemo


\Refs
\widestnumber\key{A}

\ref 
\key \cit0{1} 
\by Alekseevky, Dmitri; Kriegl, Andreas; Losik, Mark; Michor; Peter W. 
\paper Choosing roots of polynomials smoothly 
\jour Israel J. Math.
\vol 105 
\yr 1998 
\pages 203-233
\finalinfo arXiv: math.CA/9801026
\endref 

\ref
\key \cit0{2}
\by Bronshtein, M\. D\.  
\paper Smoothness of polynomials depending on parameters
\jour Sib. Mat. Zh. 
\vol 20
\yr 1979
\pages 493-501
\lang Russian
\transl\nofrills English transl. in 
\jour Siberian Math. J
\vol 20 
\yr 1980 
\pages 347-352 
\endref

\ref
\key \cit0{3}
\by Kriegl, Andreas; Michor, Peter W. 
\paper Differentiable perturbation of unbounded operators
\jour Math. Ann. (to appear)
\finalinfo arXiv: math.FA/0204060
\endref

\ref
\key \cit0{4}
\by Ohya, Y; Tarama, S.
\paper Le probl\`eme de Cauchy \`a caracteristiques multiples dans 
la classe de Gevrey (coefficients H\"olderiens en $t$)
\inbook Hyperbolic equations and related topics, Proc. Taniguchi Int. Symp., 
Katata and Kyoto/Jap. 1984 
\publ Academic Press 
\publaddr Boston
\pages 273-306 
\yr 1986
\endref

\ref
\key \cit0{5}
\by Wakabayashi, S.
\paper Remarks on hyperbolic polynomials
\jour Tsukuba J. Math.
\vol 10
\yr 1986
\pages 17--28
\endref

\endRefs
\enddocument